\documentclass[article]{siamart220329}

\usepackage{lipsum}
\usepackage{bm}
\usepackage{amsfonts}
\usepackage{graphicx}
\usepackage{epstopdf}
\usepackage{algorithmic}
\usepackage{subcaption}
\usepackage{multirow}
\usepackage{booktabs} 
\ifpdf%
  \DeclareGraphicsExtensions{.eps,.pdf,.png,.jpg}
\else
  \DeclareGraphicsExtensions{.eps}
\fi
\usepackage{amsopn}

\usepackage{booktabs}
\newcommand{\TheTitle}{%
  Generalized Optimal AMG Convergence Theory for Stokes Equations Using Smooth Aggregation and Vanka Relaxation Strategies
}

\newcommand{\TheShortTitle}{%
   Optimal AMG Convergence Theory for Stokes Equations
}

\newcommand{\TheName}{%
  Ahsan Ali
}

\newcommand{\TheAddress}{%
  Department of Mathematics and Statistics, University of New Mexico
  (\email{ahsan@unm.edu}).
}

\newcommand{\TheFunding}{%
This work was funded by NSF grant DMS-2110917.
}

\newcommand{\TheCollaborators}{%
  James J. Brannick, Karsten Kahl, Oliver A. Krzysik, Jacob B. Schroder, Ben S. Southworth, and Alexey Voronin
}

\author{\TheName\thanks{\TheAddress}}
\title{{\TheTitle}\thanks{\TheFunding}}
\headers{\TheShortTitle}{\TheName}
\ifpdf%
\hypersetup{%
  pdftitle={\TheTitle},
  pdfauthor={\TheName}
}
\fi

\begin{document}

\maketitle
\begin{center}
In collaboration with:
  {\TheCollaborators}
\end{center}
\vspace{1cm}

\begin{abstract}
  This paper discusses our recent generalized optimal algebraic multigrid (AMG) convergence theory applied to the steady-state Stokes equations discretized using Taylor-Hood elements ($\pmb{ \mathbb{P}}_2/\mathbb{P}_{1}$). The generalized theory is founded on matrix-induced orthogonality of the left and right eigenvectors of a generalized eigenvalue problem involving the system matrix and relaxation operator. This framework establishes a rigorous lower bound on the spectral radius of the two-grid error-propagation operator, enabling precise predictions of the convergence rate for symmetric indefinite problems, such as those arising from saddle-point systems. We apply this theory to recently developed monolithic smooth aggregation AMG (SA-AMG) solver for Stokes, constructed using evolution-based strength of connection, standard aggregation, and smoothed prolongation. The performance of these solvers is evaluated using additive and multiplicative Vanka relaxation strategies. Additive Vanka relaxation constructs patches algebraically on each level, resulting in a nonsymmetric relaxation operator due to the partition of unity being applied on one side of the block-diagonal matrix. Although symmetry can be restored by eliminating the partition of unity, this compromises convergence. Alternatively, multiplicative Vanka relaxation updates velocity and pressure sequentially within each patch, propagating updates multiplicatively across the domain and effectively addressing velocity-pressure coupling, ensuring a symmetric relaxation. We demonstrate that the generalized optimal AMG theory consistently provides accurate lower bounds on the convergence rate for SA-AMG applied to Stokes equations. These findings suggest potential avenues for further enhancement in AMG solver design for saddle-point systems.
\end{abstract}


\section{Introduction}\label{sec:intro} Differential equations and discrete problems with a saddle-point structure are prevalent in numerous scientific and engineering applications. A well-known example is the incompressible Stokes equations, where the incompressibility constraint results in saddle-point problems both at the continuum level and for commonly used discretizations~\cite{elman2014finite,benzi2005numerical}. Stable finite-element discretizations of such systems, such as Taylor-Hood elements ($\pmb{ \mathbb{P}}_k/\mathbb{P}_{k-1}$), lead to a $2\times2$ block structure in the resulting linear system, where the lower-right diagonal block is a zero matrix. This indefiniteness, combined with the coupling between velocity and pressure unknowns, presents significant challenges for the direct application of standard geometric or algebraic multigrid (AMG) methods.

Several successful preconditioners and solvers for saddle-point systems have been developed, bootstrap algorithm~\cite{brannick2014bootstrap}, approximate block-factorization techniques~\cite{elman2014finite, adler2017preconditioning} or monolithic geometric multigrid approaches~\cite{brandt1979multigrid, vanka1986block, linden2005multigrid, braess1997efficient, schoberl2003schwarz}. Monolithic multigrid methods, in particular, have demonstrated effectiveness by employing simultaneous coarsening of velocity and pressure variables and tailored relaxation schemes such as Vanka~\cite{vanka1986block} and Braess-Sarazin~\cite{braess1997efficient}. However, extending these geometric methods to robust and efficient AMG frameworks for a wide range of Stokes discretizations remains a challenge. While prior work on AMG for Stokes systems exists~\cite{wabro2004coupled, metsch2013algebraic, janka2008smoothed, prokopenko2017algebraic}, achieving optimal convergence for such symmetric indefinite systems is still a developing area~\cite{notay2016new, bacq2022new}.

In this study, we leverage the recently introduced generalized optimal algebraic multigrid (AMG) convergence theory \cite{ali2024generalized} to assess the performance of the state-of-the-art monolithic smooth aggregation AMG (SA-AMG) solver for the Stokes equations \cite{voronin2023monolithic}. This innovative monolithic AMG solver not only achieves but often exceeds the efficiency of existing preconditioners, demonstrating robust performance across a wide range of problem types. However, the solver's effectiveness is highly sensitive to the selection of smoother damping parameters, which are frequently determined through computationally intensive brute-force searches tailored to specific problem classes and multigrid configurations.

The optimal AMG convergence theory \cite{ali2024generalized} establishes rigorous lower bounds on the spectral radius of the two-grid error-propagation operator. This framework relies on matrix-induced orthogonality of the eigenvectors associated with the system matrix and relaxation operator. Using this theoretical foundation, we analyze the convergence properties of the monolithic SA-AMG solver \cite{voronin2023monolithic} when employing both additive and multiplicative Vanka relaxation methods. These approaches handle velocity-pressure coupling differently: additive Vanka utilizes algebraically defined patches with nonsymmetric relaxation, while multiplicative Vanka sequentially updates components to enforce symmetric relaxation. Our findings demonstrate the predictive strength of the generalized optimal AMG theory and its potential to uncover deeper insights into the convergence behavior of AMG solvers for indefinite systems.

The structure of this paper is as follows: \Cref{sec:disc} introduces the governing equations and finite element discretizations for the Stokes equations. \Cref{sec:mono-sa-amg} reviews the recently developed monolithic SA-AMG preconditioner and its features tailored for the Stokes problem. \Cref{sec:theory_frame} discusses the generalization of the optimal AMG theory framework designed for both nonsymmetric and symmetric indefinite problems. \Cref{sec: analysis} presents the overall two-grid AMG convergence analysis of the monolithic SA-AMG solver. Finally, the paper concludes with \Cref{sec:conclusion}.

\section{Governing Equations and Discretization}\label{sec:disc}

The steady-state Stokes equations, describing the flow of incompressible viscous fluids, are considered on a bounded Lipschitz domain $\Omega \subset \mathbb{R}^2$. The equations are given by:  
\begin{subequations}\label{eq:stokes-eq}
\begin{alignat}{2}
-\nabla^2 \mathbf{u} + \nabla p &= \mathbf{f}, \quad && \text{in } \Omega, \\
-\nabla \cdot \mathbf{u} &= 0, \quad && \text{in } \Omega, \\
\mathbf{u} &= \mathbf{g}_{\text{D}}, \quad && \text{on } \Gamma_{\text{D}}, \\
\frac{\partial \mathbf{u}}{\partial \mathbf{n}} - \mathbf{n} p &= \mathbf{g}_{\text{N}}, \quad && \text{on } \Gamma_{\text{N}},
\end{alignat}
\end{subequations}  
where $\mathbf{u}$ is the velocity, $p$ is the pressure, $\mathbf{f}$ is a forcing term, and $\mathbf{g}_{\text{D}}$, $\mathbf{g}_{\text{N}}$ are the boundary conditions on the Dirichlet and Neumann segments of the boundary, respectively. Here, $\mathbf{n}$ denotes the outward unit normal vector.  

To discretize the equations, we use Taylor-Hood finite elements, with piecewise quadratic basis functions for velocity and piecewise linear basis functions for pressure. The weak formulation seeks $\mathbf{u} \in \bm{\mathcal{V}}_h$ and $p \in \mathcal{Q}_h$ such that  
\begin{subequations}\label{eq:stokes-disc}
\begin{alignat}{2}
a(\mathbf{u}, \mathbf{v}) + b(p, \mathbf{v}) &= F(\mathbf{v}), \quad && \forall \mathbf{v} \in \bm{\mathcal{V}}_h, \\
b(q, \mathbf{u}) &= 0, \quad && \forall q \in \mathcal{Q}_h,
\end{alignat}
\end{subequations}  
where the bilinear forms $a(\cdot, \cdot)$ and $b(\cdot, \cdot)$ are defined as  
\begin{equation}
    a(\mathbf{u}, \mathbf{v}) = \int_\Omega \nabla \mathbf{u} : \nabla \mathbf{v}, \quad b(p, \mathbf{v}) = -\int_\Omega p \nabla \cdot \mathbf{v},  
\end{equation}
and the linear form $F(\mathbf{v})$ includes contributions from external forces and Neumann boundary conditions:  
\begin{equation}
    F(\mathbf{v}) = \int_\Omega \mathbf{f} \cdot \mathbf{v} + \int_{\Gamma_{\text{N}}} \mathbf{g}_{\text{N}} \cdot \mathbf{v}.
\end{equation}

On a structured grid in two dimensions, the discretization leads to a saddle-point linear system:  
\begin{equation}\label{eq:saddle}
K\mathbf{x} =
\begin{bmatrix}
A & B^T \\
B & 0
\end{bmatrix}
\begin{bmatrix}
\mathbf{u} \\ p
\end{bmatrix}
=
\begin{bmatrix}
\mathbf{f} \\ 0
\end{bmatrix}=\mathbf{b},
\end{equation}  
where $A$ represents the discrete vector Laplacian, and $B$ corresponds to the discrete divergence operator. For Taylor-Hood elements, the number of velocity unknowns exceeds the number of pressure unknowns. This block structure, derived from the finite-element discretization, forms the basis for our analysis of the AMG solver's convergence properties.

We use the standard Taylor-Hood (TH) discretization, denoted by {$\pmb{ \mathbb{P}}_k/\mathbb{P}_{k-1}$}, which involves continuous piecewise polynomials of degree $k$ for the velocity space $\bm{\mathcal{V}}_h$ and continuous piecewise polynomials of degree $k-1$ for the pressure space $\mathcal{Q}_h$. The TH element pair is known to be inf-sup stable for $k \geq 2$ on any triangular mesh $\Omega_h$ of the domain $\Omega$. In this study, we focus on the case where $k = 2$, which corresponds to the lowest-order pair of TH elements, with quadratic velocity components and linear pressure elements.

\section{Monolithic SA-AMG for Stokes}\label{sec:mono-sa-amg} 

In this section, we describe a state-of-the-art monolithic SA-AMG solver~\cite{voronin2023monolithic} designed for solving the Stokes system $K\mathbf{x} = \mathbf{b}$. This approach incorporates several advanced techniques, including: (i) the construction of algebraic multigrid hierarchies for different components of the system, such as the vector Laplacian governing the velocity variables (denoted by $A$) and the \textit{div-grad}-type operator for the pressure; (ii) a tailored algebraic coarsening strategy for the \((\mathbf{u}, p)\) degrees of freedom (DoFs) that replicates the behavior of efficient geometric multigrid methods for Stokes problems; and (iii) an algebraic Vanka-type relaxation scheme for saddle-point systems, with an optimized relaxation damping parameter. 

While many of these techniques are highly effective for elliptic PDE systems like the vector Laplacian, they cannot be directly applied to saddle-point systems. Standard SA-AMG methods, such as those available in libraries like PyAMG~\cite{BeOlScSo2023}, generally rely on assumptions that are invalid for saddle-point problems. For example, these methods assume that damped Jacobi relaxation is effective for smoothing errors and converging when applied to the system matrix. However, this is not the case for saddle-point systems, as the zero diagonal block makes such relaxation schemes undefined.

To address these challenges, \cite{voronin2023monolithic} employ a block-structured approach for grid transfer operators. By defining separate transfer operators for the velocity and pressure fields, the algorithm accommodates the unique characteristics of each field. For the velocity variables, this uses elliptic nature of the vector Laplacian, discretized with \(\pmb{\mathbb{P}}_1\) elements, to reliably construct interpolation operators using standard SA-AMG techniques. For the pressure field, where the system matrix contains a zero diagonal block, an auxiliary operator is required. Therefore the SA-AMG setup opt for the pressure Laplacian, \(A_p\), discretized using canonical \(\mathbb{P}_1\) pressure basis functions, as it avoids the complexity of algebraically truncating the \(BB^T\) operator.

The interpolation operators for each field are constructed independently using SA-AMG and combined into a block-diagonal monolithic operator at each level. The coarse-grid operator is computed recursively using these monolithic interpolation operators, forming a grid hierarchy that extends to the coarsest shared level. A critical aspect of this approach is maintaining stability across coarse grids to ensure convergence. This requires consistent coarsening rates between velocity and pressure fields. Traditional strength-of-connection (SoC) measures often demand significant parameter tuning to achieve this balance. Instead, the algorithm uses evolution SoC measure, which integrates local representations of algebraically smooth error, to ensure uniform coarsening rates without extensive parameter adjustments.

The standard two-level error propagation operator for a
multigrid cycle is given by
\begin{equation}
  {E_{TG}} = (I - \omega M^{-1} K)^{\nu_2} (I - P
    (RKP)^{-1} R K) ( I -\omega M^{-1}
    K)^{\nu_1}  , \label{eq:TG-operator-full}
\end{equation}
where
$K$ is the  system matrix,
$M$ represents the relaxation operator
applied to $K$ with damping parameter $\omega$,
$P$ is the interpolation 
and $R= P^T$ is the restriction operator. A key component of the monolithic SA-AMG method is the coupled relaxation scheme, characterized by the weighted error propagation operator $I - \omega M^{-1}K$. Vanka relaxation \cite{vanka1986block}, a well-established method for saddle-point problems, employs an overlapping domain-decomposition approach to iteratively refine approximations to the global solution.  While geometric multigrid methods naturally construct such patches using topological arguments, the monolithic SA-AMG solver constructs these patches algebraically at each level of the hierarchy. Consequently, the solver relies on the additive form of Vanka relaxation. More details of this novel SA-AMG algorithm and its features can be found in \cite{voronin2023monolithic}.

\section{Generalized AMG Convergence Theory Review}\label{sec:theory_frame} This section discusses the most recent development in optimal AMG convergence theory \cite{ali2024generalized}, which targets nonsymmetric and indefinite problems. This theory utilizes a specific matrix-induced orthogonality of the left and right eigenvectors of a generalized eigenvalue problem involving the system matrix and the relaxation operator. By applying this generalization of the optimal convergence theory, it is possible to derive a measure of the spectral radius of the two-grid error transfer operator that is mathematically analogous to the derivation in the symmetric positive definite setting for optimal interpolation, which relies on energy norms. Specifically, Theorem 5.1 and Corollary 5.2 in \cite{ali2024generalized} can be summarized as follows:

\textit{Let $K\in\mathbb{C}^{n\times n}$ be a general non-singular matrix. Given $n_{c}$ and non-singular smoother $M$ such that $M^{-1}K$ is diagonalizable, Consider the left and right generalized eigenvectors, $V_l=\left[\mathbf{\check{v}_{i}}\right]_{i=1}^{n}$ and $V_r=\left[\mathbf{\hat{v}_{i}}\right]_{i=1}^{n}$, of the matrix pencil $(K,M)$, respectively, where the corresponding eigenvalues $\left\{\lambda_{i}\right\}_{i=1}^{n}$ are ordered such that $|1-\lambda_{1}|\geq |1-\lambda_{2}| \geq \cdots  \geq|1-\lambda_{n}|$.  Define the optimal interpolation $P_{\sharp}$ and restriction $R_{\sharp}$ to satisfy
\begin{align*}
\operatorname{range}(P_{\sharp}) &= \operatorname{range}\left(
\begin{pmatrix} \mathbf{\hat{v}}_{1} & \mathbf{\hat{v}}_{2} & \cdots & \mathbf{\hat{v}}_{n_{c}} \end{pmatrix}\right),\\
    \operatorname{range}(R_{\sharp}^{*}) &= \operatorname{range}\left(
\begin{pmatrix} \mathbf{\check{v}}_{1} & \mathbf{\check{v}}_{2} & \cdots & \mathbf{\check{v}}_{n_{c}} \end{pmatrix}\right),
\end{align*}
Then, the spectral radius of the two-grid error transfer operator \begin{equation}\label{eq:theory_identity}
     E_{TG}(P,R)=(I-P(RKP)^{-1}RK)(I-M^{-1}K),
\end{equation} is given by
\begin{equation}\label{identity_nonsymmetric}
    \rho\left(E_{TG}(P_{\sharp},R_{\sharp})\right)= |1 - \lambda_{n_{c}+1}|.
\end{equation}}

The above theory result can be further specialized to real-valued Hermitian (potentially indefinite) operators, such as those arising from saddle point systems, which will allow a natural Galerkin $R = P^T$ set of optimal transfer operators. This is Corollary 5.3 in \cite{ali2024generalized}, which reads as follows:

\textit{Let $K\in\mathbb{R}^{n\times n}$ be a Hermitian matrix. Given $n_{c}$ and a non-singular Hermitian relaxation operator $M\in\mathbb{R}^{n\times n}$ such that $M^{-1}K$ is diagonalizable, consider the generalized EVP
 \begin{equation}\label{gen_prob_hermitian}
K\mathbf{v}_{i}=\lambda_{i}M\mathbf{v}_{i}, 
 \end{equation}
where $i=1,2,\dots,n$; and the corresponding eigenvalues $\left\{\lambda_{i}\right\}_{i=1}^{n}$ are ordered such that $|1-\lambda_{1}|\geq |1-\lambda_{2}| \geq \cdots  \geq|1-\lambda_{n}|$. Define the optimal interpolation operator $P_{\sharp}$ to satisfy
\begin{equation*}
\operatorname{range}(P_{\sharp}) = \operatorname{range}\left(
\begin{pmatrix} \mathbf{v}_{1} & \mathbf{v}_{2} & \cdots & \mathbf{v}_{n_{c}} \end{pmatrix}\right),
\end{equation*}
where we assume that when complex conjugate eigenpairs are considered, if one complex conjugate eigenvector is placed in $\operatorname{range}(P_{\sharp})$, then the other must be at well. Then letting $R_{\sharp}=P_{\sharp}^{*}$, the spectral radius of the two-grid error transfer operator is given by
\begin{equation}\label{identity}
\rho\left(E_{TG}(P_{\sharp})\right)= |1 - \lambda_{n_{c}+1}|.
\end{equation}}
For more details about this theory framework, see \cite{ali2024generalized}.

\section{AMG Convergence Analysis}\label{sec: analysis} In this section, we present a direct comparison of the performance of the standalone SA-AMG solver in terms of convergence with our theoretical framework of optimal AMG convergence theory. Two different operators must be addressed to achieve this comparison: the Stokes operator 
$K$ and the relaxation operator $M$.
First, the Stokes operator $K$ is symmetric and indefinite, but due to its construction, the additive Vanka relaxation operator $M$ is nonsymmetric. This lack of symmetry arises from the fact that the patches are assembled algebraically. Since the Dirichlet degrees of freedom (DoFs) are not algebraically connected to any pressure DoFs, they do not appear in the patches, and their contributions are missing from the relaxation matrix. Due to the nonsymmetry of the additive Vanka relaxation operator $M$, we rely on the nonsymmetric optimal AMG theory to compare the geometric convergence factor obtained from a $V(1,0)$ cycle.

\begin{table}[H]
    \centering
    \caption{Performance of SA-AMG with Different Parameters}
    \label{tab:sa_amg_results}
    \begin{tabular}{ccccc}
        \toprule
        $V$-Cycle & $\omega$ & DoFs & Iterations & Residual \\
        \midrule
        $V(1,0)$ & 0.30 & 1095 & $100^{*}$ & 2.88e-09 \\
                 &       & 1891 & $100^{*}$ & 1.74e-09 \\
                 & 0.41  & 1095 & 85        & 9.06e-11 \\
                 &       & 1891 & 83        & 7.94e-11 \\
                 & 0.52  & 1095 & 65        & 8.58e-11 \\
                 &       & 1891 & 63        & 8.50e-11 \\
                 & 0.62  & 1095 & 52        & 7.66e-11 \\
                 &       & 1891 & 50        & 8.60e-11 \\
                 & 0.73  & 1095 & 48        & 9.97e-11 \\
                 &       & 1891 & 55        & 8.73e-11 \\
                 & 0.84  & 1095 & $100^{*}$ & 3.57e-07 \\
                 &       & 1891 & $100^{*}$ & 3.32e-05 \\
                 
        \midrule
        $V(2,0)$ & 0.30  & 1095 & 60        & 9.85e-11 \\
                 &       & 1891 & 59        & 7.48e-11 \\
                 & 0.41  & 1095 & 43        & 8.66e-11 \\
                 &       & 1891 & 42        & 7.28e-11 \\
                 & 0.52  & 1095 & 46        & 6.49e-11 \\
                 &       & 1891 & $100^{*}$ & 1.62e-10 \\
                 & 0.62  & 1095 & $100^{*}$ & 9.71e-08 \\
                 &       & 1891 & $100^{*}$ & 7.92e+01 \\
                 & 0.73  & 1095 & $100^{*}$ & 7.54e+02 \\
                 &       & 1891 & $100^{*}$ & 9.75e+09 \\
                 & 0.84  & 1095 & $100^{*}$ & 1.14e+10 \\
                 &       & 1891 & $100^{*}$ & 9.73e+15 \\
                 
        \midrule
        $V(2,2)$ & 0.30  & 1095 & 32        & 9.39e-11 \\
                 &       & 1891 & 62        & 7.51e-11 \\
                 & 0.41  & 1095 & 60        & 8.18e-11 \\
                 &       & 1891 & $100^{*}$ & 1.64e-06 \\
                 & 0.52  & 1095 & 65        & 8.09e-11 \\
                 &       & 1891 & $100^{*}$ & 7.41e-07 \\
                 & 0.62  & 1095 & 38        & 9.90e-11 \\
                 &       & 1891 & 80        & 9.45e-11 \\
                 & 0.73  & 1095 & 20        & 5.29e-11 \\
                 &       & 1891 & 37        & 6.03e-11 \\
                 & 0.84  & 1095 & 23        & 5.67e-11 \\
                 &       & 1891 & 22        & 3.37e-11 \\
                 
        \midrule
        $V(4,0)$ & 0.30  & 1095 & 34        & 7.14e-11 \\
                 &       & 1891 & 63        & 7.34e-11 \\
                 & 0.41  & 1095 & 59        & 8.79e-11 \\
                 &       & 1891 & $100^{*}$ & 1.90e-06 \\
                 & 0.52  & 1095 & 64        & 8.36e-11 \\
                 &       & 1891 & $100^{*}$ & 7.14e-07 \\
                 & 0.62  & 1095 & 37        & 6.57e-11 \\
                 &       & 1891 & 79        & 8.88e-11 \\
                 & 0.73  & 1095 & 20        & 8.14e-11 \\
                 &       & 1891 & 37        & 8.76e-11 \\
                 & 0.84  & 1095 & 21        & 9.05e-11 \\
                 &       & 1891 & 22        & 4.01e-11 \\
                 
        \bottomrule
    \end{tabular}
    \label{tab:sa-amg-vcycle}
\end{table}

As with many multigrid relaxation schemes, Vanka relaxation is known to be sensitive to the choice of the damping parameter, $\omega$. The monolithic SA-AMG setup described in \cite{voronin2023monolithic} adopts a practical approach, where relaxation on each level is accelerated using the FGMRES algorithm. This is applied in two distinct scenarios: (i) In the outer iteration, where AMG serves as the preconditioner. For each FGMRES iteration, one AMG $V$-cycle is applied, (ii) FGMRES is also utilized to encapsulate Vanka relaxation, effectively acting as a makeshift Chebyshev iteration. This approach eliminates the need to establish Chebyshev bounds for each level within the multigrid hierarchy. The solver fixes the relaxation to include two inner FGMRES iterations, with one Vanka relaxation sweep per iteration, effectively setting $\nu_1 = \nu_2 = 2$ in~\eqref{eq:TG-operator-full}.  

It has been determined that FGMRES acceleration is only advantageous when performing at least two FGMRES iterations (i.e., a polynomial of order 2 or higher). For a single relaxation sweep, it is more effective to use a fixed damping parameter  and deactivate FGMRES acceleration. Therefore, a more thorough parameter search for $\omega$ is necessary that ensures $h$-independent convergence. Table \ref{tab:sa-amg-vcycle} demonstrates these results where maximum number of iterations is set to 100.

The results highlight the impact of relaxation damping parameters ($\omega$) and $V$-cycle configurations on SA-AMG convergence. Lower $\omega$ values (\(\omega \leq 0.62\)) generally lead to faster convergence, with residuals reducing to machine precision in fewer iterations. However, higher $\omega$ values ($\omega \geq 0.73$) significantly degrade convergence, causing larger residuals or divergence. $V$-cycle configurations also affect performance, with increased pre- and post-relaxations improving convergence for smaller $\omega$, as seen in the $V(2,2)$ and $V(4,0)$ cycles. However, for larger $\omega$, even these configurations fail to stabilize, resulting in poor convergence . This underscores the sensitivity of SA-AMG solver to parameter tuning.

\begin{figure}[h!]
     \centering
     \begin{subfigure}[b]{0.49\textwidth}
         \centering
\includegraphics[width=\textwidth]{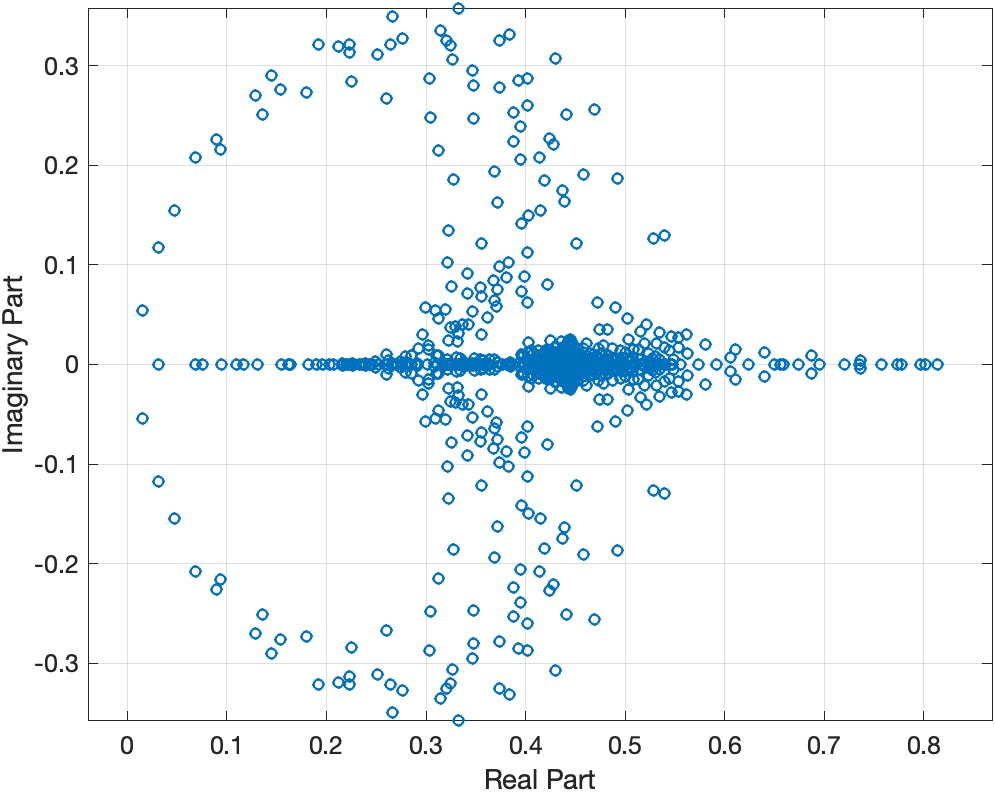}
         \caption{$\omega=0.41$}
         \label{fig:eig_1}
     \end{subfigure}
     \hfill
     \begin{subfigure}[b]{0.49\textwidth}
         \centering
         \includegraphics[width=\textwidth]{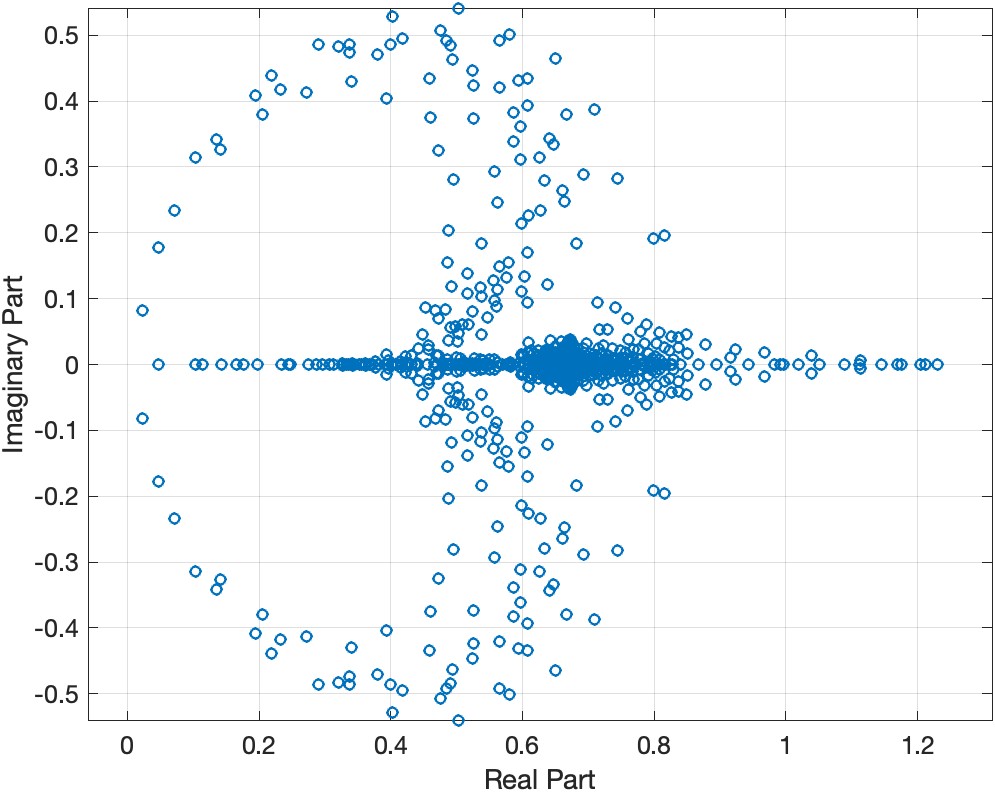}
         \caption{$\omega=0.62$}
         \label{fig:eig_2}
     \end{subfigure}
     \hfill
     \begin{subfigure}[b]{0.49\textwidth}
         \centering
         \includegraphics[width=\textwidth]{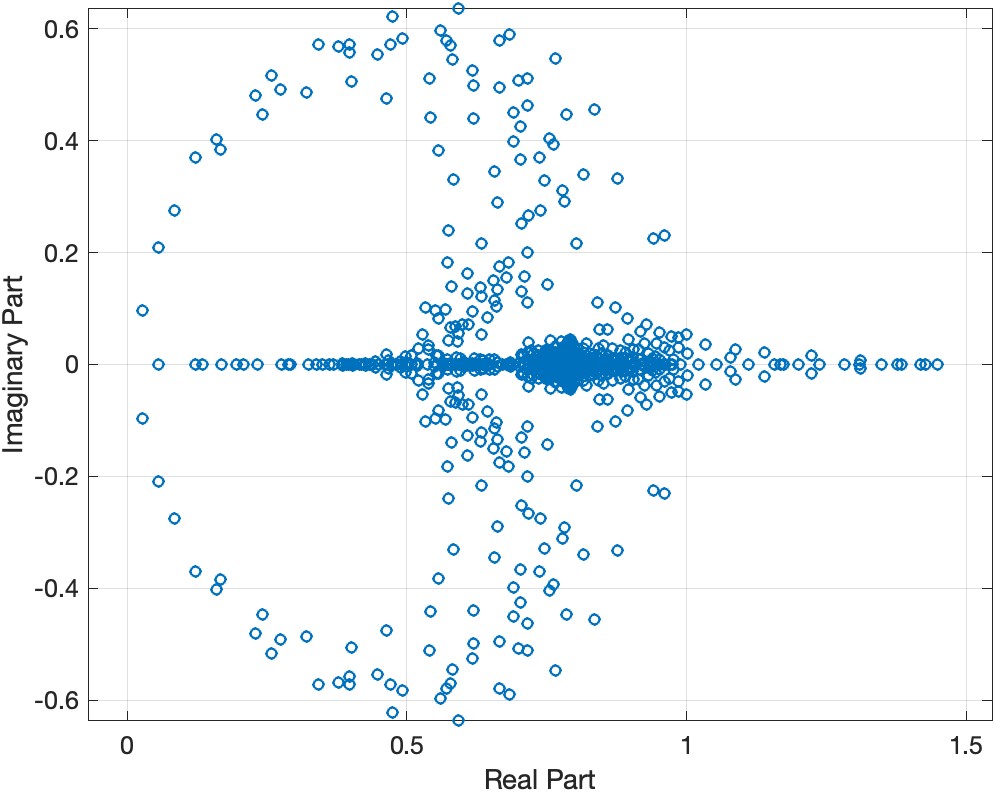}
         \caption{$\omega=0.73$ }
         \label{fig:eig_3}
         \end{subfigure}
         \hfill
     \begin{subfigure}[b]{0.49\textwidth}
         \centering
         \includegraphics[width=\textwidth]{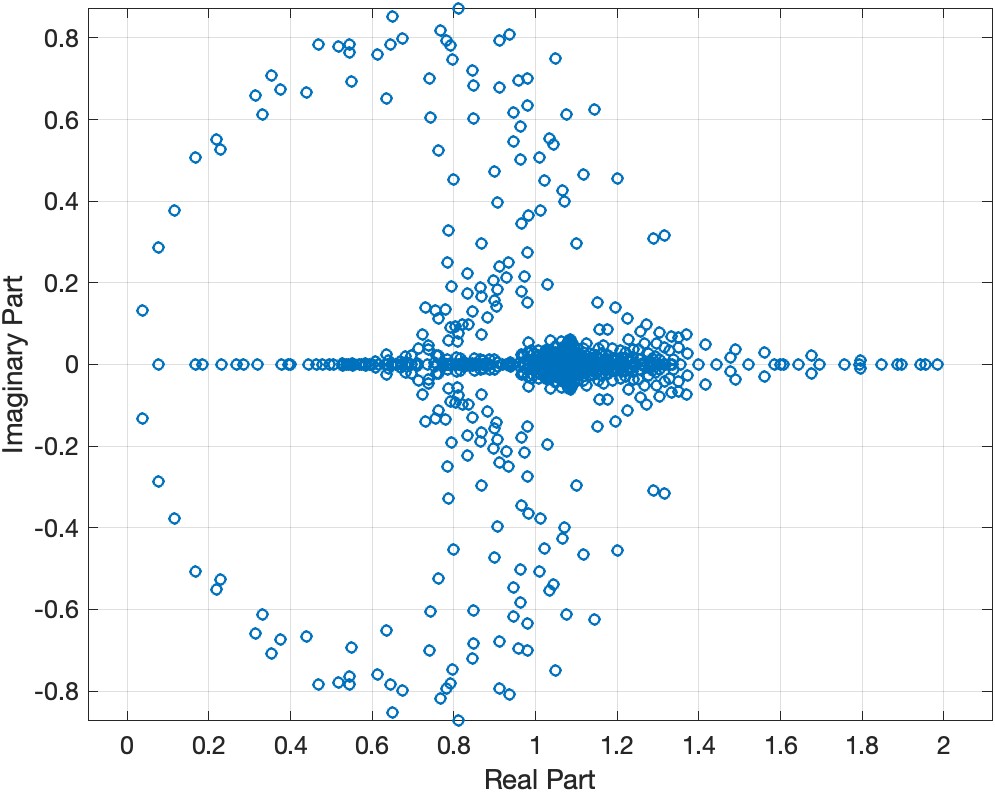}
         \caption{$\omega=1$  }
         \label{fig:eig_4}
         \end{subfigure}
    
    \caption{Spectrum of complex eigenvalues for different damping parameters.}
    \label{fig:spectrum}
\end{figure}

In order to fully understand how damping parameter effects spectrum of eigenvalues of the generalized eigenvalue problem $K\mathbf{v}=\lambda M \mathbf{v}$, we illustrate the spectrum of the eigenvalues for fixed Stokes operator $K$ of size $1095 \times 1095$, Additive Vanka relaxation operator $M$, with different $\omega$. Figure \ref{fig:spectrum} implies that widely scattered eigenvalues can cause the iterative method to stagnate, making it difficult to achieve convergence within a reasonable number of iterations.

\begin{figure}[h!]
  \centering
  \includegraphics[width=0.7\textwidth]{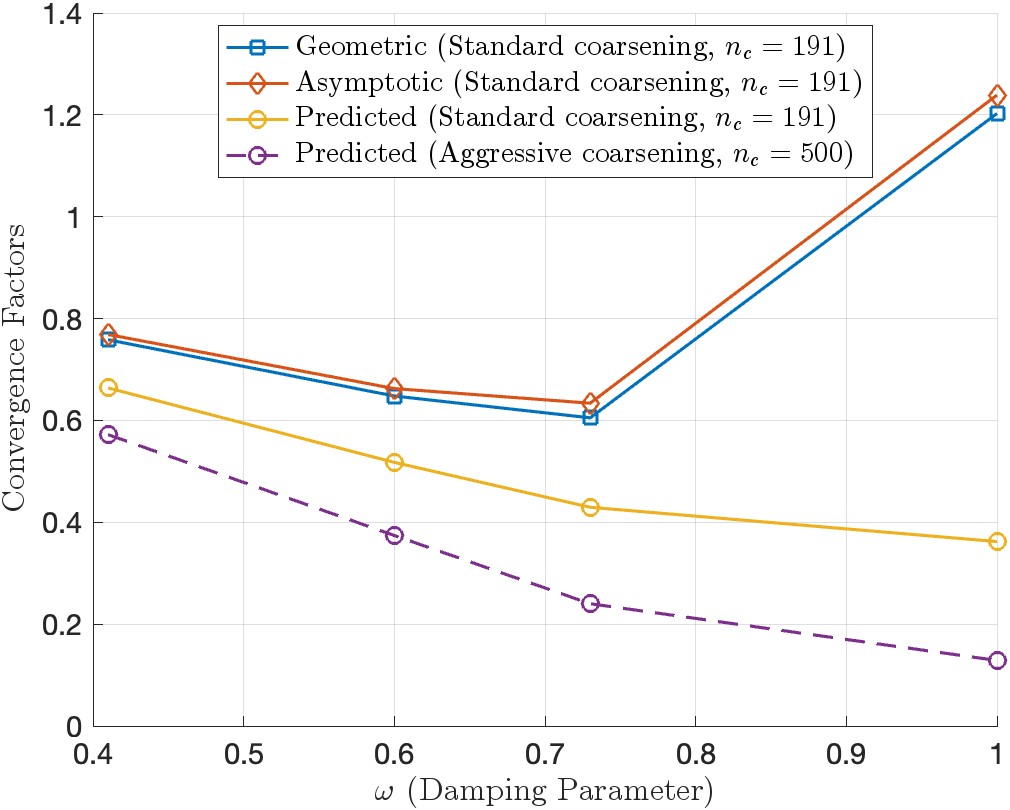}
  \caption{Plot of the observed geometric and asymptotic convergence factors with standard coarsening automatically chosen through the monolithic SA-AMG setup. The observed convergence factors are compared with the predicted convergence determined using the optimal AMG convergence theory identity $|1-\lambda_{n_{c}+1}|$
, which consistently provides a lower bound and serves as a theoretical framework to improve convergence.}
  \label{fig:convergence comparsion}
\end{figure}

Figure \ref{fig:convergence comparsion} demonstrates the comparison of geometric and asymptotic convergence factors for the two-grid monolithic SA-AMG  applied to a Stokes problem. The focus is on a $V(1,0)$ cycle setup, as our theoretical framework for optimal interpolation is developed under the assumption of one pre-smoothing step and no post-smoothing. For this test, the Stokes operator $K$ has a size of $1095 \times 1095$, and the diagonal block interpolation operator $P$ has a size of $1095 \times 191$. Consequently, for standard coarsening, we set  $n_{c} = 191$ as the number of coarse points. The relaxation method used is the additive Vanka relaxation operator \(M\), which is nonsymmetric due to the partition of unity. Unlike setups that encapsulate Vanka relaxation within FGMRES, this study directly evaluates the performance of the  SA-AMG as a standalone solver. Specifically, the geometric and asymptotic convergence factors are compared against the theoretical predictions derived from the identity \(|1 - \lambda_{n_c + 1}|\). 

Key observations indicate that the geometric and asymptotic convergence factors align closely across different damping parameters \(\omega\), demonstrating the robustness of the AMG setup. The predicted convergence factors, based on the theoretical identity \(|1 - \lambda_{n_c + 1}|\), consistently exhibit lower values, underestimating the observed convergence factors and suggesting potential areas for improvement. Using \(n_c = 500\), which corresponds to aggressive coarsening, the theoretical identity predicts smaller asymptotic convergence factors, indicating faster convergence. These results confirm that the generalized optimal AMG theory provides an accurate framework for predicting convergence rates in AMG methods, particularly those designed for saddle-point systems. Additionally, the theory's ability to capture convergence behavior under both standard and aggressive coarsening scenarios highlights its potential to inform and improve AMG solver design for saddle-point problems.

\begin{figure}[h!]
  \centering
  \includegraphics[width=0.75\textwidth]{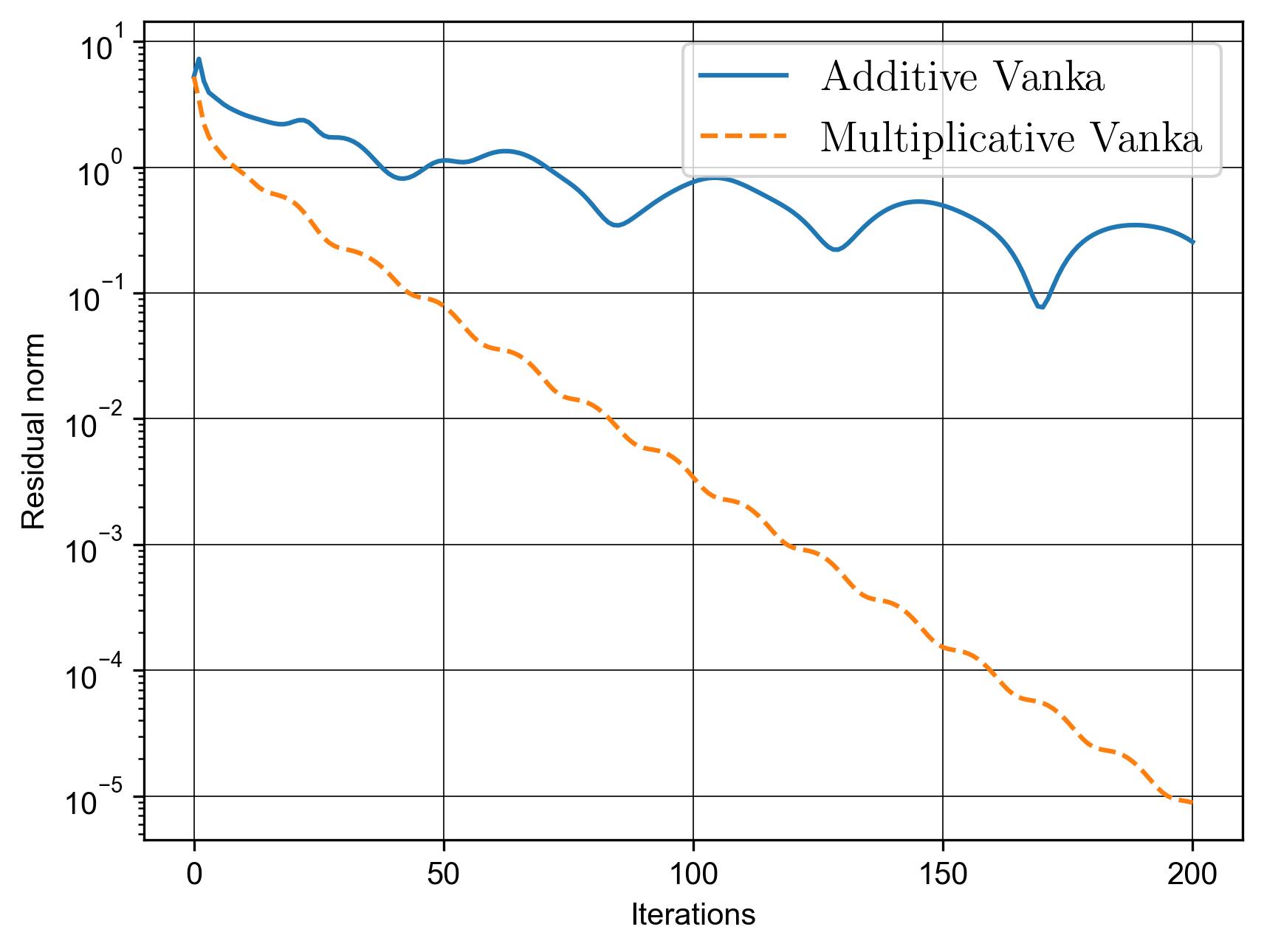}
  \caption{ Convergence behavior of two Vanka relaxation methods, Additive Vanka and Multiplicative Vanka, applied to solve the linear system $K\mathbf{x}=\mathbf{b}$ arising from the Stokes problem.}
  \label{fig:vanka_compare}
\end{figure}

In this SA-AMG setting, Multiplicative Vanka can show  more potential in achieving faster convergence. Figure \ref{fig:vanka_compare} compares the convergence behavior of the Additive and Multiplicative Vanka relaxation methods for solving a linear system. The plot shows the residual norm as a function of iterations for both methods. The Multiplicative Vanka method, indicated by the dashed line, exhibits a faster reduction in the residual norm compared to the Additive Vanka method, which is represented by the solid line. The steeper initial decline in residuals for the Multiplicative Vanka method suggests that it more efficiently reduces the error per iteration, while the Additive Vanka method demonstrates a slower but steady convergence. Both methods ultimately converge to similar residual levels, with Multiplicative Vanka achieving this in fewer iterations.

\section{Conclusion and Future Work}\label{sec:conclusion} 
In this work, we show that the generalized optimal AMG theory provides an accurate framework for predicting convergence rates in AMG methods designed for saddle point systems. The ability of the theory to capture the convergence behavior under both standard and aggressive coarsening scenarios highlights its potential for advancing AMG solver design for saddle-point problems. Future work will explore leveraging these findings to refine coarsening and smoothing strategies, refine existing state-of-the-art solvers, and develop even more efficient AMG solvers for indefinite systems. In particular, we will focus on compatible relaxation \cite{brannick2010compatible} to efficiently determine the smoother parameter and develop an efficient generalized reduction-based AMG solver \cite{ali2024constrained} for indefinite problems.
\bibliographystyle{siamplain}
\bibliography{refs}

\end{document}